\documentclass[a4paper,12pt]{article}

\usepackage{amsmath,amssymb,algorithmic}
\usepackage{latexsym}
\usepackage[latin1]{inputenc}
\usepackage[T1]{fontenc}
\usepackage{textcomp}
\usepackage{times}

\usepackage[pdftex]{graphicx}
\usepackage[frenchb]{babel}
\usepackage{babelbib}
\usepackage{hyperref}
\hypersetup{
colorlinks=true, 
citecolor=blue, %
pdfmenubar=false,
pdftoolbar=false,
bookmarksopen=true,
breaklinks=true, 
urlcolor=blue, 
linkcolor=blue, 
}

\author{Philippe Besse\thanks{Université de Toulouse -- INSA, Institut de Mathématiques, UMR CNRS 5219} \and Aurélien Garivier\thanks{Université de Toulouse -- UT3, Institut de Mathématiques, UMR CNRS 5219} \and  Jean-Michel Loubes\thanks{Université de Toulouse -- UT3, Institut de Mathématiques, UMR CNRS 5219}}
\title{\emph{Big Data Analytics} -- Retour vers le Futur~-~3~- De~Statisticien à \emph{Data Scientist}}

\begin{document}
\DeclareGraphicsExtensions{.pdf,.jpg}
\maketitle

\begin{quote}
{\bf Résumé}: L'évolution rapide des systèmes d'information gérant des données de plus en plus volumineuses a causé de profonds changements de paradigme dans le travail de statisticien, devenant successivement prospecteur de données, bio-informaticien et maintenant \emph{data scientist}. Sans souci d'exhaustivité et après avoir illustré ces mutations successives, cet article présente brièvement les nouvelles questions de recherche qui émergent rapidement en Statistique, et plus généralement en Mathématiques, afin d'intégrer les caractéristiques: volume, variété et vitesse, des grandes masses de données.

{\bf Mots-clefs}: Statistique; Fouille de Données; Données Biologiques à Haut Débit; Grande Dimension; Bioinformatique; Apprentissage Statistique; Grandes Masses de Données.

{\bf Abstract}: The rapid evolution of information systems managing more and more voluminous data has caused profound paradigm shifts in the job of statistician, becoming successively data miner, bioinformatician and now data scientist. Without the concern of completeness and after having illustrated these successive mutations, this article briefly introduces the new research issues that quickly arise in Statistics, and more generally in Mathematics, in order to integrate the characteristics: volume, variety and velocity, of big data.

\end{quote}

\section{Introduction: du kO au PO}
Il était une fois la Statistique : une question (i.e. biologique), associée à une hypothèse expérimentalement réfutable, une expérience planifiée où $p$ (au plus quatre) variable et facteurs sont observées sur $n$ (environ trente) individus, un modèle linéaire supposé vrai, un test, une décision, une réponse. Ce cadre rassurant de la Statistique du milieu du siècle dernier a subi de profonds changements de paradigme, conséquences des principaux sauts technologiques ; en schématisant, à chaque fois que le volume des données à analyser est multiplié par un facteur mille. 

Sans viser l'exhaustivité, cet article s'attache à analyser ou plutôt illustrer les changements de paradigmes qui accompagnent les formes successives que prend le métier de statisticien, au gré des développements brutaux des outils informatiques et, notamment, de ceux des systèmes d'information. Comment le statisticien est-il passé du traitement de quelques kilo-octets à celui de méga, giga, téra et maintenant péta-octets ? Quels sont les principaux changements méthodologiques qui ont accompagné ces mutations, et quelles en sont leurs conséquences sur les thèmes de recherche ?

Trois exemples sont développés dans une première section pour illustrer ce propos :
\begin{itemize}
\item	fouille de données (\emph{data mining}) appliquée à la gestion de la relation client (GRC) de la fin des années 90 ;
\item	apprentissage statistique de données omiques la décennie suivante ; 
\item	avènement récent et très médiatisé du \emph{big data}.
\end{itemize}
Nous terminons en pointant de nouvelles pistes de recherche qui s'ouvrent, et en précisant les compétences nécessaires à un statisticien pour devenir un \emph{data scientist}.

Ce court article ne vise pas une présentation générale des très vastes thèmes concernés par les données massives et leur valorisation.  Il tente d'expliciter, à partir de l'expérience et d'exemples, le point de vue de mathématiciens / statisticiens sur les outils de la composante \emph{Analytics} souvent absente des présentations. Il s'agit ainsi de comprendre et mesurer l'impact de ces nouvelles technologies  tant sur les enseignements universitaires que sur certains aspects de la recherche à venir. 

``Le consensus aujourd'hui est de définir le \emph{data scientist} à l'intersection de trois domaines d'expertise: Informatique, Statistique et Mathématiques, et Connaissances métier'' (\href{http://abiteboul.com/DOCS/14.DataScience.pdf}{Abitboul et al; 2013}). Cet article propose, du point de vue du statisticien, une incitation au débat entre ces trois domaines. 

\section{Trois changements de paradigme}
\subsection{Exploration et fouille de MO de données}
Le premier changement de paradigme est lié à l'origine des données, elles deviennent préalables à l'analyse avec le \emph{data mining} (Fayyad et al. 1996). L'expérience n'est plus planifiée en relation avec l'objectif, la question posée, les données sont issues d'entrepôts, de bases et systèmes relationnels, stockées pour d'autres raisons notamment comptables. Dans un grand nombre de cas, les données massives ne sont pas utilisées pour produire des connaissances : l'objectif est d'apprendre directement des règles de décision ou de diagnostic efficientes. Moyennant quelques connaissances du langage de requête SQL, le statisticien, devenu \emph{prospecteur}, conçoit les méthodes qui classifient ou segmentent des bases de clientèle, recherchent des cibles, qu'elles soient à privilégier et solliciter, à travers des scores d'appétence ou d'attrition, ou bien à éviter, en estimant des risques de défaut de paiement, de faillite. 

La fouille de données n'a pas suscité de développements méthodologiques particulièrement originaux. Pilotés à l'aide d'une interface graphique plus ou moins ergonomique, des assemblages logiciels dédiés, libres ou commerciaux, intègrent langage de requête dans les bases de données, exploration statistique et réduction de dimension (analyse de données en France ou \emph{multivariate analysis}), de classification non supervisée ou \emph{clustering}, de modélisation statistique classique (régression gaussienne et logistique, analyse discriminante) et plus algorithmique (arbres de décision binaires, réseaux de neurones). Cet ensemble est complété pour analyser le « panier de la ménagère » par la recherche de règles d'association. 

Le \emph{data mining} s'est ainsi beaucoup développé à l'initiative des entreprises du secteur tertiaire (banque assurance, VPC, téléphonie...) pour des objectifs de marketing quantitatif, très peu dans les entreprises industrielles ou la Statistique restait une obligation réglementaire (cf. l'industrie pharmaceutique) plus qu'une nécessité économique. 

\subsection{Apprentissage en ultra haute dimension de GO de données}
Le deuxième changement de paradigme est illustré par l'avènement de l'imagerie et de celui des biotechnologies. Il a fallu 15 ans et environ 2,7 milliards de dollars pour finir de séquencer un génome humain en 2003. La même opération s'exécute en 2014 en quelques heures pour quelques milliers d'euros. Le prix n'est plus un goulot d'étranglement, ce sont les capacités d'analyse qui le deviennent. Aux 3,2 milliards de paires de nucléotides, s'ajoutent les dizaines voire centaines de milliers d'informations omiques relatives à l'expression des gènes (ARN), à celles des protéines (spectrométrie de masse) ou encore des métabolites (résonance magnétique nucléaire) présents dans un échantillon biologique. Ces données sont stockées, annotées, dans des bases dédiées et souvent en accès public. Du point de vue du biologiste, les tâches et fonctions d'analyse sont confondues mais concernent tant le bioinformaticien qui aligne et compare des séquences, les dénombre, reconstitue des génomes, que le biostatisticien qui cherche à extraire de ces données des informations significatives. 

Le biologiste, qui observe ainsi des dizaines de milliers ($p$) de variables sur quelques dizaines ($n$) d'échantillons, voit croître le volume d'informations, tandis que le statisticien ne voit croître que l'indétermination du système où le nombre de variables est considérablement plus important que celui des individus. Cette situation extrême a suscité,  justifié, beaucoup de développements méthodologiques concernant les corrections de tests multiples (FDR de Benjamini et Hochberg ; 1995), la recherche de modèles parcimonieux ou \emph{sparse}, l'utilisation d'algorithmes d'apprentissage (SVM) de complexité liée à $n$ pas à $p$ (Guyon et al. 2002) ou insensibles (\emph{random forest}) au sur-apprentissage (Genuer et al. 2010) couplés à une procédure de sélection de variables.

Certes, les masses de données générées et stockées sont considérables mais, une fois filtrées, normalisées, pré-traitées, par le bioinformaticien, le biostatisticien, qui peut être la même personne, n'analyse qu'un tableau de taille globale relativement classique bien que très largement disproportionné par le nombre de colonnes ou variables. C'est	ainsi que la plupart des méthodes statistiques usuelles, de la régression aux modèles graphiques en passant par la PLS (i.e. Lê Cao et al. 2008 ; 2011) ont été revisitées par l'adjonction d'un terme de pénalisation dite \emph{Lasso} (Tibshirani, 1996), ou en norme $L_1$, afin de forcer l'annulation de nombreux paramètres et donc arriver à une sélection drastique des variables conjointement à l'estimation du modèle. L'objectif est d'obtenir des résultats simplement interprétables ou encore, avec un objectif prévisionnel, c'est la recherche de biomarqueurs, protéines, gènes, susceptibles, par exemple, d'affiner un diagnostic précoce. Au cours de nombreux projets, les auteurs ont apprécié l'importance de la ``connaissance métier'' relevée par Abitboul et al. (2013) qui rend indispensable une étroite collaboration avec les biologistes dont les expériences ultimes valident une classification, une sélection de gènes, un réseau et donc la méthodologie développée et mise en place.

Le principal souci du statisticien devenu \emph{bioinformaticien} est le contrôle du sur-ajustement en apprentissage statistique (Hastie et al. 2009) afin d'éviter les artefacts de l'échantillon. Ceci se traduit par la recherche d'un compromis classique entre biais et variance, entre erreur d'approximation et d'estimation, pour optimiser la complexité du modèle et ainsi minimiser le risque ou erreur de prévision ou de généralisation.

\subsection{Science des grandes masses de TO de données}
Dans les applications en production industrielle, le e-commerce, la géo-localisation... c'est le nombre $n$ d'individus et la complexité des données qui explosent. Les bases se structurent en nuages (\emph{cloud}), les moyens de calculs se groupent (\emph{cluster}), mais la puissance ne suffit pas à la voracité (\emph{greed}) des algorithmes. Encouragé par la loi de Moore, de plus en plus est attendu des méthodes de traitement automatique des données, et on constate que le nombre de variables observées sur chaque individu croît lui aussi ($p$ grandit avec $n$). Dans le savant équilibre qui conduit aux règles de décision les plus pertinentes, un troisième terme d'erreur vient s'ajouter à l'approximation et à l'estimation : l'erreur  d'optimisation, induite par la limitation du temps de calcul ou celle du volume / flux de données considéré. 

En relation avec ces nouvelles structures de données réparties ou fragmentées (\emph{sharding}) en n\oe uds (cf.\href{http://hadoop.apache.org/}{ \emph{Hadoop}}) le statisticien devenu \emph{data scientist} revisite ses classiques pour se focaliser sur des outils et méthodes conduisant à des exécutions \emph{échelonnables} (\emph{scalable}) et compatibles (\href{http://mahout.apache.org/}{\emph{Mahout}} , \href{https://github.com/RevolutionAnalytics/RHadoop/wiki}{\emph{R-Hadoop}} ) avec les nouveaux systèmes d'information. S'il doit développer des algorithmes ainsi distribués, il (ré)apprend la syntaxe des langages fonctionnels (Lisp, Scheme, Caml) pour utiliser \href{http://www.scala-lang.org/}{\emph{Scala}} ou \href{http://clojure.org/}{\emph{Clojure}}.

Le traitement de données volumineuses favorise donc un retour en grâce des langages (fonctionnels) et des méthodes robustes : $k$-means (Mac Queen, 1967), classifieur bayésien naïf (Maron; 1961) facilement distribuables (MapReduce) sur des centaines voire milliers de n\oe uds associant un processeur et un ensemble de données. On observe en outre un vif regain d'intérêt pour les problèmes d'allocation dynamique de ressources appelés « bandit manchot », étudiés depuis longtemps (Thomson, 1933) dans les essais cliniques, pour optimiser les systèmes de recommandation (par exemple dans le e-commerce). Ces problèmes sont très liés au « A/B testing », une technique connue en marketing visant à améliorer une réponse (par exemple un taux de clic) en comparant un groupe de contrôle à un groupe test. Il est aujourd'hui courant de pouvoir effectuer et analyser un nombre très important de tels tests, pour lesquels il est à la fois essentiel, de maximiser l'information reçue pour chaque observation et de contrôler la significativité globale des résultats.

\section{Les pistes de recherche}
Le \href{http://www.mckinsey.com/insights/business_technology/big_data_the_next_frontier_for_innovation}{rapport} de McKinsey (2011)  a popularisé la caractérisation du \emph{big data} par trois « V » (volume, variété, vélocité). Ceux-ci ne font pas que faire revisiter le passé, ils propulsent le mathématicien / statisticien / \emph{data scientist} vers de nouveaux horizons de recherche mais aussi dans un écosystème, une jungle aux imbrications matérielles, logicielles, économiques... excessivement complexes. La \href{http://mattturck.com/2012/06/29/a-chart-of-the-big-data-ecosystem/}{figure 1} due à  Turk et Zillis (2012) est un instantané sûrement partiel illustrant bien la complexité du paysage des acteurs des grandes masses de données. Il mêle éditeurs commerciaux de grands logiciels « ravalant » et interfaçant des outils bien rodés, start-ups très discrètes sur leurs algorithmes et fondations de code source libre largement diffusés ; tous développant des infrastructures, des outils d'analyse, des applications ou produisant des données. Premiers concernés, les acteurs du monde informatique, se mobilisent rapidement dans le sillage des leaders du marché (Google, Yahoo, Amazon, Facebook) du e-commerce.

De leur côté, les acteurs, mathématiciens, statisticiens, évoluent comme c'est souvent le cas avec d'autres échelles de temps. D'une part ils développent des cadres théoriques les plus fondamentaux pour une meilleure compréhension des problèmes étudiés. D'autre part, ils proposent de nouvelles méthodologies et algorithmes dont ils évaluent les propriétés pour cerner leurs domaines de validité : tantôt ils précèdent l'application, tantôt ils théorisent a posteriori des  algorithmes et méthodes pour en étudier les propriétés. 

\begin{figure}
\centerline{\includegraphics[width=11cm]{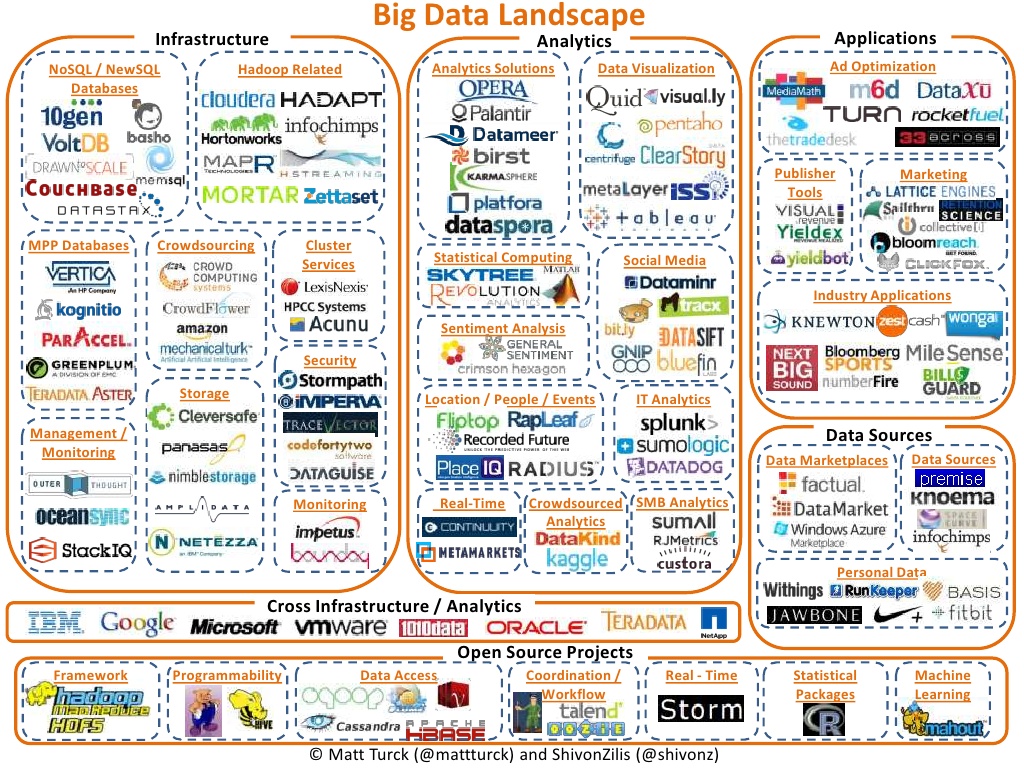}}
\caption{\it Carte de l'échosystème du Big Data en juin 2012.}\label{jungle}
\end{figure}

\subsection{Apprentissage et statistique}
Le statisticien académique du siècle dernier a largement développé les études théoriques afin de qualifier les méthodes utilisées en caractérisant les lois des paramètres, majorant finement les risques de prévision et cherchant les tests les plus puissants. La dynamique de recherche change avec les développements des algorithmes et méthodes à l'interface entre Apprentissage Machine et Statistique.

 Symptomatiquement, la « méta » librairie caret (Kuhn; 2008) de R propose l'accès et la comparaison de près de 140 méthodes de modélisation statistique ou d'apprentissage en régression ou classification supervisée, elles-mêmes réparties dans plusieurs autres librairies de R. Certains modèles (gaussien, binomial...) sont étudiées depuis très longtemps, alors que d'autres méthodes comme les forêts aléatoires (\emph{random forest} de Breiman ; 2001) sont toujours l'objet de recherches actives afin d'en préciser ou justifier les bonnes propriétés (i.e. Biau ; 2012). Pour beaucoup de ces méthodes, en particulier celles qui s'appuient sur le principe de minimisation empirique du risque, des résultats sont obtenus à partir d'inégalités de concentration probabilistes (i.e. Birgé et Massart ; 1993), elles-mêmes parfois issues d'analyses géométriques (i.e. Ledoux et Talagrant ;1991). Bien entendu, avant que puissent être prouvées les propriétés d'un algorithme d'apprentissage, ses performances pratiques sont évaluées sur des données synthétiques ou réelles de bases publiques comme celle de l'UCI (Bache et Lichman ; 2013). Reconnaissons par ailleurs que pour beaucoup de ces comparaisons\footnote{Voir à ce propos les scénarios d'analyse de plusieurs jeux de données publics proposés sur le site \href{http://wikistat/fr}{Wikistat}.}, les différences observées entre ces performances ne sont pas toujours très significatives et qu'il n'est pas pertinent, pour un exemple pratique donné, de se lancer dans une comparaison systématique de toutes ces techniques!
 
Il reste beaucoup de travail aux statisticiens et probabilistes pour étudier les propriétés des méthodes d'apprentissage en concurrence mais l'arrivée des données massives change radicalement l'approche des problèmes en soulevant de nombreuses questions et donc pistes de recherche. Sans prétendre à l'exhaustivité, en voici quelques unes triées selon les caractéristiques des trois « V ».

\subsection{Volume}
L'infrastructure dominante \emph{Hadoop} pour distribuer et archiver les données est issue directement des besoins de Google (archiver des pages, compter des occurrences de mots...) et de son système de fichier.  Faire ce choix, et donc celui de la parallélisation par \emph{MapReduce}, impacte fortement les modes d'analyse susceptibles d'intégrer le changement d'échelle en fonction de l'objectif poursuivi. Les approches en marketing peuvent utiliser directement des méthodes échelonnables (\emph{scalable}) comme celle des $k$ plus proches voisins ou faire appel à un algorithme stochastique dérivé de celui de Robins et Monro (1951) et plus généralement d'algorithmes de descente de gradient stochastique (ou non) pour pouvoir estimer un modèle bien connu comme une régression logistique.  Mais finalement, très peu de méthodes d'apprentissage satisfont directement aux contraintes d'une parallélisation par \emph{MapReduce}. Owen et al. (2012), Giacomelli (2013), le site \href{https://mahout.apache.org/}{\emph{Mahout}}, celui de \href{http://blog.cloudera.com/blog/2013/03/cloudera_ml_data_science_tools/}{\emph{cloudera ML}}, listent les rares directement interfacées avec \emph{Hadoop}: $k$-means, classifieur naïf bayésien, \emph{random forest}, régression logistique. 

D'autre part, quand ces méthodes sont utilisées, une grande partie du travail très approfondi de ces dernières années en Apprentissage Machine et Statistique semble oubliée ou tout du moins est absente des livres référencés. Les recommandations concernant la taille de l'échantillon ou sous-échantillon à considérer ne sont guère précises dans la littérature de l'environnement \emph{Hadoop}; elles sont discutées vis-à-vis des capacités de calcul et de stockage, pas en termes de qualité de prévision. C'est pourtant cette réflexion qui permet d'optimiser un modèle (sélection de variables, paramètre de pénalisation ou de complexité), choisir une méthode et utiliser des outils diagnostics évaluant qualité et pertinence d'une aide à la décision. Il serait important de comparer, par exemple, l'erreur de prévision d'une technique exécutant une optimisation stochastique échelonnable, avec celle d'une méthode simplement séquentielle estimée et surtout \emph{optimisée} à partir d'un  sondage aléatoire simple ou stratifié. 

Ceci nécessite la disponibilité d'un algorithme d'échantillonnage. Il est souvent fait référence au \emph{reservoir sampling} proposé par Vitter (1985) mais ScaSRS (\emph{scalable simple random sampling}) proposé par Xiangrui (2013) semble plus approprié. Enfin, \emph{random forest}  (Breiman; 2001) est la rare technique récente d'apprentissage facilement échelonnable et donc implémentée dans une librairie comme \emph{Mahout}. Les bonnes propriétés de cette méthode (résistance au sur apprentissage) et son algorithme d'agrégation de modèles la rendent  particulièrement adaptée.  Les modèles (arbres de décision) estimés sur des échantillons, en principe \emph{bootstrap}, le sont sur des échantillons plus  ``indépendants'' si la taille des données augmentent. Elle devient la méthode (classification ou régression) non-linéaire de référence mais ne peut généralement se substituer à d'autres lorsque le phénomène à modéliser est linéaire ou qu'une interprétation détaillée du modèle est nécessaire.

\subsection{Variété}
Le grand nombre de données étudiées entraîne de fait une très grande hétérogénéité. Sous le terme de variété se cache donc de multiples problématiques. 

Tout d'abord les données, notamment dans un environnement industriel, peuvent se présenter sous des formes variées telles que des textes, séquences, fonctions, courbes, spectres, images, graphes...  voire des combinaisons de l'ensemble. Au delà des difficultés techniques inhérentes au traitement d'une telle diversité, se pose le problème de l'utilisation de toute cette masse d'information hétéroclite.   Chaque structure de donnée soulève des questions originales spécifiques et, comparer ces différentes structures entre elles, devient également compliqué à la fois d'un point de vue pratique mais également théorique. En effet calculer ne serait-ce qu'une simple moyenne ou des distances entre les objets concernés requiert une analyse particulière pour prendre en compte leur géométrie. Il est par exemple possible de prendre en considération des approximations de la distance géodésique sur une variété riemanienne (Diméglio et al. 2013).

 D'autre part ces données sont observées avec une grande variabilité qui masque bien souvent la nature du phénomène étudié. Il importe dès lors de révéler la structure contenue dans ces observations afin d'en extraire l'information qu'elles contiennent. Ainsi, pour calculer la distance entre deux courbes, celle usuelle (espaces $L_2$) est beaucoup trop sensible à ces légères déformations ou décalages des courbes. Un recalage préalable (time warping) est nécessaire. Une autre approche prometteuse (scattering) en analyse d'images consiste à représenter celles-ci dans une base d'ondelettes possédant des propriétés d'invariance pour certains groupes de transformations (Mallat et Sifre ; 2012). Cette représentation permet alors d'identifier ou regrouper des formes ou textures d'images particulières. 
 
Le dernier exemple concerne initialement l'analyse d'images, conjointement avec des principes de parcimonie. Candès et al. (2004) ont montré que le \emph{compressed sensing} permet d'acquérir beaucoup plus rapidement, avec moins de mesures, une très bonne précision des images médicales à partir du moment où celles-ci sont structurées. La reconstruction de l'image fait ensuite appel à la résolution parcimonieuse d'un système linéaire (\emph{Dantzig selector}, Candès et Tao, 2005) par optimisation convexe en norme $L_1$. Celle-ci peut avoir bien d'autres applications notamment pour l'analyse des très grandes matrices très creuses, produites en génétique, génomique, fouille de textes, analyse d'incidents dans l'industrie... et leur complétion (Candès et Tao; 2009).

\subsection{Vélocité}
La vélocité des données est également motrice de nouvelles recherches sur les algorithmes des méthodes de décision qui deviennent nécessairement adaptatives ou séquentielles. Comment traiter au mieux le flux de données au fil de l'eau ? Le cadre classique de l'apprentissage en ligne ne s'y prête guère, qui supposait un échantillon fixé, pour lequel le statisticien avait le temps de faire tous les calculs nécessaires afin d'obtenir une règle de décision qui ne changerait pas de sitôt. 

Les approches intrinsèquement séquentielles, permettant de s'adapter progressivement au nouvelles données qui arrivent, connaissent ainsi un vif regain d'intérêt ; pour l'étape d'optimisation conduisant aux règles de décision, c'est notamment le cas de la descente de gradient stochastique (Bach et Moulines 2013). En outre, il convient de prendre en compte le vieillissement des données : si, en première approche, on peut penser à utiliser une simple fenêtre glissante ou un facteur d'escompte (Garivier et Moulines 2011), des traitements intrinsèquement séquentiels et adaptatifs sont plus convaincants. \`A titre d'exemple, on peut penser aux systèmes de recommandation, comme ceux qui interviennent dans la gestion du contenu des pages web : le succès fulgurant de \href{http://www.criteo.com/fr}{Criteo}  souligne l'importance stratégique de cette problématique. Par ailleurs, l'organisation et le choix des prix de réserves dans les systèmes d'enchères au second prix à grande vitesse (Cesa-Bianchi et al.; 2013) nécessitent également la mise en production ultra rapide et décentralisée d'algorithmes issus de modèles probabilistes avancés qui rappellent les problématiques du \emph{trading} à haute fréquence.

Ce sujet illustre tout particulièrement le fait que, dans un cadre séquentiel, les décisions prises à un instant peuvent influencer les observations futures, qui du coup ne peuvent pas être traitées comme un échantillon. Ainsi, le système de recommandation ne peut se contenter de proposer les contenus dont l'appétence estimée est la plus grande car, s'il le fait, il se trouvera vite piégé, ne servant plus qu'un nombre très restreint de contenus pas assez diversifiés -- les autres n'étant simplement pas assez proposés pour que l'intérêt que leur porte l'utilisateur puisse être perçu. 

Les modèles de bandit manchot (White ; 2012), qui intéressent depuis très longtemps les concepteurs d'essais cliniques, constituent ici un modèle de référence, pour lequel Cappé et al. (2013) proposent un contrôle optimal efficace. Ce modèle admet un certain nombre de variantes, mais il mérite encore d'être enrichi et complété pour pouvoir répondre à toute la diversité des contraintes applicatives. Récemment, il a été montré que des méthodes très simples d'inspiration bayésienne pouvaient aussi avoir un comportement optimal (Kaufmann et al.; 2012), ce qui laisse beaucoup d'espoirs quant à la possibilité de traiter des modèles complexes.

\section{Conclusion}
Ce rapide tour d'horizon montre comment, sous la pression des données, le métier de statisticien évolue selon le contexte, en prospecteur,  bioinformaticien, \emph{data scientist}. Mais, à travers ces avatars, l'objectif principal reste celui d'apporter un quatrième « V » à la définition du big data, celui de la Valorisation. Seules des aides à la décision pertinentes sont susceptibles de justifier les efforts financiers consacrés à l'acquisition et au stockage devenu massif des données. 

L'épopée décrite dans cet article insiste sur des interactions de plus en plus fortes et collaborations très intriquées entre Informatique et Mathématiques pour relever au mieux les nouveaux défis posés par l'afflux la complexité et le volume des données. A l'intérieur même des mathématiques, les contributions tendent à se diversifier. Le \emph{machine learning} s'est ainsi beaucoup appuyé sur la modélisation stochastique et statistique des données afin de construire des algorithmes fournissant des règles de décisions pertinentes et pour obtenir des garanties théoriques d'efficacité ou d'optimalité de ces méthodes en se fondant principalement sur des outils probabilistes (parfois d'inspiration géométrique).  Les grandes masses de données et de signaux suscitent un grand intérêt pour l'étude des matrices aléatoires de grande dimension, et en particulier de leur spectre asymptotique alors que les premiers résultats ont été initiés dans ce domaine par des questions de physique théorique (i.e. Johnstone ; 2006). Le passage à l'échelle du \emph{Big Data} rend indispensables et centraux de nouveaux travaux dans le domaine de l'optimisation convexe et plus généralement de l'analyse ;  la modélisation des systèmes complexes, à commencer par les grands réseaux, suscitent des rapprochements avec certains travaux de géométrie ou d'algèbre.  Le rôle du mathématicien y est important pour apporter un changement de perspective, souvent contre-intuitif mais efficace, en élevant le niveau d'abstraction ou intégrant une approche stochastique. Cela permet, par exemple, de réduire le volume des mesures : du plan d'expérience au \emph{compressed sensing}, de montrer la convergence vers un optimum global (algorithme stochastique), de construire le bon critère pour mesurer une distance adaptée (géodésique) ou introduire une pénalisation  ($L_1$) de parcimonie.

Pertinence et efficacité de l'analyse rendent nécessaires des compétences en Informatique, Mathématiques comme en Statistique. Elles ne peuvent raisonnablement être acquises qu'au niveau Master 2 ou même doctorat et des équipes pluridisciplinaires sont indispensables. Au statisticien, ou plutôt maintenant au \emph{data scientist}, d'assurer l'interface entre des compétences en systèmes d'information, langage (java, python, scala...), algorithmique, d'une part et mathématiques d'autre part. Cela est indispensable à l'identification d'informations significatives (au sens statistique), prédictives, et à forte valeur ajoutée.

Hormis le travail des fondations et des acteurs des logiciels libres, les compagnies commerciales du e-commerce et pionnières de l'analyses des grandes masses de données sont remarquablement discrètes sur leurs pratiques, les algorithmes et programmes utilisés. En conséquence la recherche ``académiques'' doit parallèlement se développer, notamment aux interfaces info/maths/métiers, pour permettre à ces approches d'être transférées  et applicables à d'autres fins et dans d'autres domaines, comme celui de la Santé Publique. C'est typiquement le cas de la factorisation de matrices non négatives (NMF), largement utilisée en e-commerce et dont une implémentation en java \emph{scalable}  est rendue disponible (Ruiqi et al; 2014) grâce à  la recherche financée sur fonds publics en Biologie.

Plus fondamentalement, la transmission de ces connaissances et compétences nouvelles est d'autant plus importante à mettre en place au regard des enjeux sociétaux. Cela concerne évidemment les nombreuses questions éthiques et juridiques sur les droits de collecte et d'utilisation des données, leur propriété, le droit à l'oubli, l'évolution à venir d'internet. Le débat auquel nous invite \href{https://webwewant.org/}{Tim Berners-Lee}\footnote{https://webwewant.org/} ne peut se tenir sans des citoyens informés sur les outils qu'ils utilisent.

\section*{Bibliographie}
Abiteboul S., Bancilhon F., Clemencon S., Saporta G. (2014) L'émergence d'une nouvelle filière de formation: ``data scientists''. Rapport technique, 11 pages, (http://abiteboul.com/DOCS/14.DataScience.pdf).

Benjamini Y., Hochberg Y. (1995). Controlling the false discovery rate: a practical and powerful approach to multiple testing.  \emph{Journal of the Royal Statistical Society}, Series B,  57-1, 289-300. 

Bach F. et Moulines E. (2013). Non-strongly-convex smooth stochastic approximation with convergence rate O(1/n). To appear in \emph{Advances in Neural Information Processing Systems} (NIPS).

Bache K. et Lichman, M. (2013).  \href{http://archive.ics.uci.edu/ml}{UCI Machine Learning Repository}. Irvine, CA: University of California, School of Information and Computer Science.

Biau, G. (2012). Analysis of a random forests model, \emph{Journal of Machine Learning Research}, Vol. 13, pp. 1063-1095.

Breiman L. (2001). Random forests. \emph{Machine Learning}, 45:5-32.

Candès E. J.,  Romberg J.  et Tao T. (2004). Robust uncertainty principles: exact signal reconstruction from highly incomplete frequency information. \emph{IEEE Trans. Inform. Theory}, 52 489-509.

Candès E. J. et Tao T. (2005). The Dantzig selector: statistical estimation when p is much larger than n. \emph{Annals of Statistics}, 35 2313-2351.

Candès E. J.  et Tao. T.  (2009). The power of convex relaxation: Near-optimal matrix completion. \emph{IEEE Trans. Inform. Theory} 56(5), 2053-2080.

Cappé O., Garivier A., Maillard O. A., Munos R. et Stoltz G. (2013). Kullback-Leibler Upper Confidence Bounds for Optimal Sequential Allocation. \emph{Annals of Statistics}, vol.41(3) pp.1 516-1 541.

Cesa-Bianchi N., Gentile C., Mansour Y. (2013). Regret Minimization for Reserve Prices in Second-Price Auctions. \emph{SODA}. 1190-1204.

Dimeglio C., Gallón S., Loubes J-M., et Maza, E. (2014). A robust algorithm for template curve estimation based on manifold embedding. \emph{Computational Statistics \& Data Analysis}, 70-C, 373-386.

Fayyad U. M., Piatetsky-Shapiro G. et Smyth P. (1996). From data mining to knowledge discovery : an overview, \emph{Advances in Knowledge Disco-very and Data Mining} (U. M. Fayyad, G. Piatetsky-Shapiro, P. Smyth et R. Uthurusamy, réds.), AAAI Press/MIT Press, 1996, p. 1-34.

Garivier, A. et Moulines, E. (2011). On Upper-Confidence Bound Policies for Non-stationary Bandit Problems \emph{Algorithmic Learning Theory} n°22 Oct. 2011 pp.174-188 

Genuer R., Poggi J.-M., Tuleau-Malot C. (2010) Variable selection using Random Forests. \emph{Pattern Recognition Letters} 31:2225-2236.

Giacomelli P. (2013). \emph{Apache Mahout Cookbook}. Packt Publishing.

Guyon I., Weston J., Barnhill S., and Vapnik V.N. (2002) Gene selection for cancer classification using support vector machines. \emph{Machine Learning}, 46(1-3):389-422.

Hastie T., Tibshirani R., Friedman J. (2009). \emph{The Elements of Statistical Learning: Data Mining, Inference, and Prediction}. Second Edition, Springer.

Johnstone I. M. (2006). High dimensional statistical inference and random matrices, \emph{Proceedings of the International Congress of Mathematicians}, Madrid, Spain.

Kaufmann E., Korda N., Munos R. (2012). Thompson Sampling: an asymptotically optimal finite-time analysis. \emph{Proceedings of the 23rd International Conference on Algorithmic Learning Theory} (ALT).

Kuhn M. (2008). Building Predictive Models in R Using the caret Package, \emph{Journal of Statistical Software}, 28-5.

Lê Cao K. A., Rossouw D., Robert-Granié C., Besse P. (2008). A sparse PLS for variable selection when integrating Omics data, \emph{Statistical Applications in Genetics and Molecular Biology}, Vol. 7 : Iss. 1, Article 35. 

Lê Cao K. A., Boistard S., Besse P. (2011) Sparse PLS Discriminant Analysis: biologically relevant feature selection and graphical displays for multiclass problems. \emph{BMC Bioinformatics}, 12:253.

MacQueen J. B. (1967). Some Methods for classification and Analysis of Multivariate Observations. \emph{Proceedings of 5th Berkeley Symposium on Mathematical Statistics and Probability 1}. University of California Press. pp. 281-297.

Mallat S. et Sifre  L. (2012). Combined Scattering for Rotation Invariant Texture Analysis, \emph{Proc. European Symposium on Artificial Neural Networks}.

Marron M. E. (1961). Automatic Indexing: An Experimental Inquiry, \emph{Journal of the ACM (JACM)}, VOl. 8 : Iss. 3, pp 404-417.

Owen S., Anil R., Dunning T. et Friedman E. (2012). \emph{Mahout in action}. Manning.

Thompson W. (1933). On the likelihood that one unknown probability exceeds another in view of the evidence of two samples, \emph{Bulletin of the American Mathematics Society}, vol. 25, pp. 285-294, 1933

Ruiqi L., Yifan Z., Jihong G., Shuigeng Z. (2014) CloudNMF: A MapReduce Implementation of Nonnegative Matrix Factorization for Large-scale Biological Datasets, Genomics, Proteomics \& Bioinformatics, Volume 12, Issue 1, Pp 48-51, 
(http://www.sciencedirect.com/science/article/pii/S1672022913000752).

Tibshirani R. (1996). Regression shrinkage and selection via the lasso. \emph{J. Royal. Statist.} Soc B., 58(1) :267-288.

Vitter J. S. (2013). Random sampling with a reservoir, \emph{Journal
ACM Transactions on Mathematical Software} (TOMS) TOMS, Volume 11 Issue 1, March 1985
Pages 37-57 

White J. M. (2012). \emph{Bandit Algorithms for Website Optimization Developing, Deploying, and Debugging}, O'Reilly Media.

Xiangrui M. (2013). Scalable Simple Random Sampling and Stratified Sampling. \emph{Proc. of the International Conference on Machine Learning}, Vol. 28, Atlanta.

\end{document}